\documentstyle[12pt]{article}

\begin{document}

\begin{center}
{\large {\bf FULL CHAINS OF TWISTS FOR ORTHOGONAL ALGEBRAS}}\\[5mm]
{\bf Petr P. Kulish} \\
St.Petersburg Department of the Steklov Mathematical  Institute, \\
191011, St.Petersburg, Russia \\
{\bf Vladimir D. Lyakhovsky} \\ Theoretical
Department, St. Petersburg State University,\\ 198904, St. Petersburg,
Russia \\ {\bf Alexander A. Stolin} \\ Department of Mathematics, 
University of Goteborg,\\ S-412 96 Goteborg, Sweden \\
\end{center}

\begin{abstract}
 We show that for some Hopf subalgebras in $U_{\cal F}(so(M))$ nontrivially 
deformed by a twist ${\cal F}$ it is possible to find the nonlinear primitive 
copies. 
This enlarges the possibilities to construct chains of twists.
For orthogonal algebra $U(so(M))$ we present a method to compose the full 
chains with carrier space as large as the Borel subalgebra $B(so(M))$. 
These chains can be used to construct the new deformed Yangians.  
\end{abstract}

\section{Introduction}

Quantizations of triangular Lie bialgebras ${\bf L}$ with antisymmetric
classical $r$-matrices $r=-r_{21}$ are defined by a twisting element $%
{\cal F}=\sum f_{\left( 1\right) }\otimes f_{\left( 2\right) }\in {\cal A}
\otimes {\cal A}$ which satisfies the twist equations \cite{DR83}:
\begin{equation}
\begin{array}{l}
\label{twist-eq}\left( {\cal F}\right) _{12}\left( \mit\Delta \otimes {\rm id}
\right) {\cal F}=\left( {\cal F}\right) _{23}\left( {\rm id}\otimes \mit\Delta
\right) {\cal F}, \\ \left( \epsilon \otimes {\rm id}\right) {\cal F}=\left(
{\rm id}\otimes \epsilon \right) {\cal F}=1.
\end{array}
\end{equation}
Explicit form of the twisting element is quite important in applications
because
it provides explicit expressions
for the quantum  ${\cal R}$-matrix ${\cal R_F}={\cal F}_{21}%
{\cal F}^{-1}$ and for the twisted coproduct $\mit\Delta _{{\cal F}}\left(
{}\right) ={\cal F}\mit\Delta _{{\cal F}}\left( {}\right) {\cal F}^{-1}$.

The first nontrivial explicitly written twisting elements ${\cal F}$ were 
given in the papers
\cite{R}, \cite{OGI}, \cite{GIA} and \cite{KLM}. These twists can be
defined on the following {\em \ carrier algebra} $L$:
$$
\begin{array}{cccc}
\left[ H,E\right] =E, &  &  & \left[ H^{\prime },E\right] =\gamma ^{\prime
}E, \\
\left[ H,A\right] =\alpha A, &  &  & \left[ H^{\prime },A\right] =\alpha
^{\prime }A, \\
\left[ H,B\right] =\beta B, &  &  & \left[ H^{\prime },B\right] =\beta
^{\prime }B, \\
\left[ E,A\right] =\left[ E,B\right] =0, &  &  & \left[ A,B\right] =E \\
\alpha +\beta =1, &  &  & \alpha ^{\prime }+\beta ^{\prime }=\gamma ^{\prime
}.
\end{array}
$$
Explicit expressions for their twisting elements are
\begin{equation}
\label{elem-twists}
\begin{array}{l}
\begin{array}{l}
\mit\Phi _{
{\cal R}}=e^{H\otimes H^{\prime }}, \\ \mit\Phi _{
{\cal J}}=e^{H\otimes \sigma }, \\ \mit\Phi _{{\cal EJ}}=
\mit\Phi _{{\cal E}}\mit\Phi _{{\cal J}}=
e^{A\otimes Be^{-\beta \sigma }}e^{H\otimes \sigma },
\end{array}
\quad
\begin{array}{l}
r_{
{\cal R}}=H\wedge H^{\prime }, \\ r_{
{\cal J}}=H\wedge E, \\ r_{{\cal EJ}}=H\wedge E+A\wedge B,
\end{array}
\\
\sigma =\ln \left( 1+E\right) .
\end{array}
\end{equation}
Here $r_{\cal R} , r_{\cal J} , r_{\cal EJ}$ are the corresponding
classical $r$-matrices.

Carrier subalgebras ${\bf L}$ can be found in any simple Lie algebra $g$
of rank greater than 1.

It was demonstrated in \cite{KLO} that these twists can be composed into {\em 
chains}. They are based on the sequences of regular injections
constructed for the initial Lie algebra
$$
g_p\subset g_{p-1}\ldots \subset g_1\subset g_0=g.
$$
To form the chain one must choose an initial root $\lambda _0$ in the root
system $\mit\Lambda \left( g\right) $, consider the set $\pi $ of its {\em 
constituent roots}%
$$
\begin{array}{l}
\pi =\left\{ \lambda ^{\prime },\lambda ^{\prime \prime }\mid \lambda
^{\prime }+\lambda ^{\prime \prime }=\lambda _0;\quad \lambda ^{\prime
}+\lambda _0,\lambda ^{\prime \prime }+\lambda _0\notin \mit\Lambda \left(
g\right) \right\}  \\
\pi =\pi ^{\prime }\cup \pi ^{\prime \prime };\qquad \pi ^{\prime }=\left\{
\lambda ^{\prime }\right\} ,\pi ^{\prime \prime }=\left\{ \lambda ^{\prime
\prime }\right\} .
\end{array}
$$
and the subset $\mit\Lambda _{\mit\Lambda _0}^{\perp }$ of roots orthogonal 
to $\lambda _0$ (the corresponding subalgebra in $g$ will be denoted by
$g_{\lambda _0}^{\perp }$ ).

It was shown that for the classical Lie algebras $g$ one can always find in $%
g_{\lambda _0}^{\perp }$ a subalgebra $g_1\subseteq g_{\lambda _0}^{\perp
}\subset g_0=g$ whose generators become primitive after the extended twist ${%
\ \mit\Phi _{{\cal EJ}}}$. Such primitivization of $g_k\subset g_{k-1}$ (called
the {\em matreshka} effect \cite{KLO}) provides the possibility to compose
chains of extended twists of the type ${\mit\Phi _{{\cal EJ}}}$,
\begin{equation}
\label{chain-ini}
\begin{array}{l}
{\cal F}_{{\cal B}_{0\prec p}}=\prod_{k=0}^p\mit\Phi _{{\cal E}_k}
\mit\Phi _{{\cal J}_k}, 
\\ \mit\Phi _{{\cal E}_k}\mit\Phi _{{\cal J}_k}=\prod_{\lambda ^{\prime }\in
\pi _k^{\prime }}\exp \left\{ E_{\lambda ^{\prime }}\otimes E_{\lambda
_0^k-\lambda ^{\prime }}e^{-\frac 12\sigma _{\lambda _0^k}}\right\} \cdot
\exp \{H_{\lambda _0^k}\otimes \sigma _{\lambda _0^k}\}\,.
\end{array}
\end{equation}

Chains of twists quantize a large variety of $r$-matrices corresponding to 
Frobenius subalgebras in simple Lie algebras \cite{STO}.

\section{Construction of a full chain of twists}

The main point in the construction of a chain is the {\em invariance} of $%
g_{k+1}$ with respect to $\mit\Phi _{{\cal E}_k{\cal J}_k}$. When these
subalgebras are proper the canonical chains have only a part of $%
B^{+}\left( g\right) $ as the twist carrier subalgebra:
\begin{equation}
\label{sequence}
\begin{array}{ccccccc}
\ldots  & \subset g_{\lambda _0^k}^{\perp } & \subset g_{\lambda
_0^{k-1}}^{\perp } & \ldots  & \subset g_{\lambda _0^1}^{\perp } & \subset
g_{\lambda _0^0}^{\perp } & \subset g \\
& \cap  & \cap  &  & \cap  & \cap  & \parallel  \\
\ldots  & \subset g_{k+1} & \subset g_k & \ldots  & \subset g_2 & \subset
g_1 & \subset g_0
\end{array}
\end{equation}

We would like to demonstrate that
the effect of primitivization is {\em universal}
and extends to the whole subalgebra $g_{\lambda _0^k}^{\perp }$. It was
shown in \cite{KL} that the invariance of a subalgebra in $g_{\lambda
_0^k}^{\perp }$ is only one of the forms of the {\em primitivization} . In
general this is the existence (in the twisted Hopf algebra $U_{{\cal E}_k%
{\cal J}_k}\left( g_{\lambda _0^k}^{\perp }\right) $ ) of a primitive
subspace $V_G^{k+1}$ with the algebraic structure isomorphic to $g_{\lambda
_0^k}^{\perp }$. On this subspace the subalgebra  $g_{\lambda _0^k}^{\perp }$
is realized nonlinearly so $V_G^{k+1}$ is called {\em deformed carrier
space} \cite{KL}.

In this context the situation with the twists for $U\left( sl(N)\right) $ is
degenerate: the subalgebra $\left( sl(N)\right) _{\lambda _0^k}^{\perp }$
coincides with $\left( sl(N)\right) _{k+1}$ , i.e. $V_G^{k+1}=V_{\left(
sl(N)\right) _{\lambda _0^k}^{\perp }}$ .

In the case of $U\left( so(M)\right) $ the situation is different. Let the
root system $\mit\Lambda \left( so(M)\right) $ be
$$\left\{ \pm e_i\pm e_j\mid
i,j=1,2,\ldots M/2;i\neq j\right\}
$$
for even $M$ and
$$\left\{ \pm e_i\pm
e_j;\pm e_k\mid i,j,k=1,2,\ldots \left( M-1\right) /2;i\neq j\right\}
$$
for odd $M$. Take $e_1+e_2$ as the initial root. Here the subalgebras $%
g_{\lambda _0^{k-1}}^{\perp }$ and $g_k$ in (\ref{sequence}) are related as
follows,
$$
g_{\lambda _0^{k-1}}^{\perp }=g_k\oplus so^{\left( k\right)
}(3)=so\left( M-4 k \right) \oplus so^{\left( k\right) }(3).
$$
Consider the invariants of the vector fundamental representations of $%
g_{k+1}=so(M-4\left( k+1\right) )$ acting on $g_k$ :
\begin{equation}
\label{invodd}
\begin{array}{l}
I_{2N+1}^a=\frac 12E_a^2+\sum_{l=3}^N\left( E_{a+l}E_{a-l}\right) , \\
I_{2N+1}^{a\otimes b}=E_a\otimes E_b+\sum_{l=3}^N\left( E_{a+l}\otimes
E_{b-l}+E_{a-l}\otimes E_{b+l}\right) ,
\end{array}
\end{equation}
\begin{equation}
\label{inveven}
\begin{array}{l}
I_{2N}^a=\sum_{l=3}^N\left( E_{a+l}E_{a-l}\right) , \\
I_{2N}^{a\otimes b}=\sum_{l=3}^N\left( E_{a+l}\otimes E_{b-l}+E_{a-l}\otimes
E_{b+l}\right) ,
\end{array}
\end{equation}
The $so^{\left( k\right) }(3)$ summands are non-trivially deformed 
by $\mit\Phi _{{\cal E}_{k-1}{\cal J}_{k-1}}$:
$$
\begin{array}{cc}
\mit\Delta _{{\cal E}_{k-1}{\cal J}_{k-1}}\left( E_{1-2}^k\right) = &
E_{1-2}^k\otimes 1+1\otimes E_{1-2}^k+\left( 1\otimes e^{-\frac 12\sigma
_{1+2}^{k-1}}\right) I_{M-4k}^{1\otimes 1} \\
& +I_{M-4k}^1\otimes \left( e^{-\sigma _{1+2}^{k-1}}-1\right) , \\
\mit\Delta _{{\cal E}_{k-1}{\cal J}_{k-1}}\left( E_{2-1}^k\right) = &
E_{2-1}^k\otimes 1+1\otimes E_{2-1}^k+\left( e^{\sigma
_{1+2}^{k-1}}-1\right) \otimes I_{M-4k}^2e^{-\sigma _{1+2}^{k-1}} \\
& +\left( 1\otimes e^{-\frac 12\sigma _{1+2}^{k-1}}\right)
I_{M-4k}^{2\otimes 2}.
\end{array}
$$
According to the main principle formulated above (despite the deformed
co-structure of $V_{g_{\lambda _0^{k-1}}^{\perp }}$) the primitivization is
realized on its isomorphic image $V_G^{k+1}$ contained
in $U_{{\cal E}_{k-1}{\cal J}_{k-1}}\left(
g_{\lambda _0^{k-1}}^{\perp }\right) $. To find this deformed carrier
subspace $V_G^{k+1}$ it is sufficient to inspect the coproducts of
invariants (\ref{invodd}) and (\ref{inveven}),
$$
\begin{array}{l}
\mit\Delta _{
{\cal E}_k{\cal J}_k}\left( I_{M-4k}^1\right) = \\ =I_{M-4k}^1\otimes
e^{-\sigma _{1+2}^k}+1\otimes I_{M-4k}^1+I_{M-4k}^{1\otimes 1}\left(
1\otimes e^{-\frac 12\sigma _{1+2}^k}\right) ,
\end{array}
$$
$$
\begin{array}{l}
\mit\Delta _{
{\cal E}_k{\cal J}_k}\left( I_{M-4k}^2e^{-\sigma _{1+2}^k}\right) = \\
=I_{M-4k}^2e^{-\sigma _{1+2}^k}\otimes 1+e^{\sigma _{1+2}^k}\otimes
I_{M-4k}^2e^{-\sigma _{1+2}^k}+I_{M-4k}^{2\otimes 2}\left( 1\otimes
e^{-\frac 12\sigma _{1+2}^k}\right).
\end{array}
$$

Now one can construct the following nonlinear primitive generators
$$
\begin{array}{l}
G_{1-2}^{k+1}=E_{1-2}^k-I_{M-4k}^1, \\
G_{2-1}^{k+1}=E_{2-1}^k-I_{M-4k}^2e^{-\sigma _{1+2}}, \\
H_{1-2}^{k-1},
\end{array}
\
\begin{array}{l}
\mit\Delta _{
{\cal E}_k{\cal J}_k}\left( G_{1-2}^{k+1}\right) =G_{1-2}^{k+1}\otimes
1+1\otimes G_{1-2}^{k+1}, \\ \mit\Delta _{
{\cal E}_k{\cal J}_k}\left( G_{2-1}^{k+1}\right) =G_{2-1}^{k+1}\otimes
1+1\otimes G_{2-1}^{k+1}, \\ \mit\Delta _{{\cal E}_k{\cal J}_k}\left(
H_{1-2}^{k-1}\right) =H_{1-2}^{k-1}\otimes 1+1\otimes H_{1-2}^{k-1}.
\end{array}
$$

The subspace spanned by $\left\{
H_{1-2}^k,G_{1-2}^{k+1},G_{-1+2}^{k+1}\right\} $ forms the algebra $%
so_G^{\left( k+1\right) }(3)\approx $ $so^{\left( k+1\right) }\left(
3\right) $ :
$$
\begin{array}{l}
\left[ H_{1-2}^k,G_{1-2}^{k+1}\right] =G_{1-2}^{k+1}, \\
\left[ H_{1-2}^k,G_{2-1}^{k+1}\right] =-G_{2-1}^{k+1}, \\
\left[ G_{1-2}^{k+1},G_{2-1}^{k+1}\right] =2H_{1-2}^k.
\end{array}
$$
Therefore we obtain the deformed primitive space
$$
V_G^{k+1}\left( g_{\lambda _0^k}^{\perp }\right) =V\left( g_{k+1}\right)
\oplus V\left( so_G^{\left( k+1\right) }(3)\right) ,
$$
that can be considered as a carrier for the twists (\ref{elem-twists}). The
next extended Jordanian twist in the chain (that is defined on $g_{k+1}$)
does not touch the space $V\left( so_G^{\left( k+1\right) }(3)\right) $.
Consequently after all the steps of the chain we will still have a
primitive subalgebra%
$$
{\cal D}=\sum_{k=0}^p\,^{\oplus }so_G^{\left( k+1\right) }(3)
$$
defined on the sum of deformed spaces $V\left( so_G^{\left( k+1\right)
}(3)\right) $.

Thus in the twisted Hopf algebra $U_{{\cal B}_{0\prec p}}\left( so\left(
M\right) \right) $ one can perform further twist deformations with the
carrier subalgebra in ${\cal D}$ . The most interesting among them are the
Jordanian twists defined by
$$
\mit\Phi _{{\cal J}_k}^G=\exp \left( H_{1-2}^k\otimes \sigma _{_G}^k\right)
\quad {\rm with\qquad }\sigma _{_G}^k\equiv \ln \left(
1+G_{1-2}^{k+1}\right)
$$
This means that in the general expression for the twisting 
element ${\cal F}_{{\cal B}_{0\prec p}}$ one can insert 
in the appropriate $k\geq 0$ places
the Jordanian twisting factors defined on the deformed carrier spaces,
i.e. to perform a substitution%
$$
\begin{array}{l}
\mit\Phi _{
{\cal E}_k}\mit\Phi _{{\cal J}_k}\Rightarrow 
\mit\Phi _{{\cal J}_k}^G\mit\Phi _{{\cal E}_k}
\mit\Phi _{{\cal J}_k}\equiv \mit\Phi _{{\cal G}_k} \\ \exp \left\{
I_{M-4k}^{1\otimes 2}\left( 1\otimes e^{-\frac 12\sigma _{1+2}^k}\right)
\right\} \cdot \exp \{H_{1+2}^k\otimes \sigma _{1+2}^k\}\Rightarrow  \\
\exp \left( H_{1-2}^k\otimes \sigma _{_G}^k\right) \cdot \exp \left\{
I_{M-4k}^{1\otimes 2}\left( 1\otimes e^{-\frac 12\sigma _{1+2}^k}\right)
\right\} \cdot \exp \left( H_{1+2}^k\otimes \sigma _{1+2}^k\right)
\end{array}
$$
This gives {\em the full chain} in the following form
\begin{equation}
\label{fullchain}
\begin{array}{l}
{\cal F}_{{\cal G}_{0\prec p}}=\prod_{k=p}^0\mit\Phi _{{\cal G}_k}= \\
\prod_{k=p}^0\left( \exp \left( H_{1-2}^k\otimes \sigma _{_G}^k\right) \cdot
\exp \left\{ I_{M-4k}^{1\otimes 2}\left( 1\otimes e^{-\frac 12\sigma
_{1+2}^k}\right) \right\} \cdot \exp \{H_{1+2}^k\otimes \sigma
_{1+2}^k\}\right) .
\end{array}
\end{equation}
Obviously the additional twistings by $\mit\Phi _{{\cal J}_k}^G$ cannot be
performed before the deformation of
the corresponding spaces $V_G^{k+1}$ by the
extended Jordanian twists $\mit\Phi _{{\cal E}_k}\mit\Phi _{{\cal J}_k}$.

\section{Applications}

The previous result means that we have constructed explicit quantizations
$$
{\cal R}_{{\cal G}_{0\prec p}}=\left( {\cal F}_{{\cal G}_{0\prec p}}\right)
_{21}\left( {\cal F}_{{\cal G}_{0\prec p}}\right) ^{-1}
$$
of the following set of classical $r$-matrices:%
$$
r_{{\cal G}_{0\prec p}}=\sum_{k=0}^p\eta _k\left( H_{1+2}^k\land
E_{1+2}^k+\xi _kH_{1-2}^k\land E_{1-2}^k+I_{M-4k}^{1\land 2}\right)
$$
Here all the parameters are independent.

The dimensions of the nilpotent subalgebras $N^{+}\left( so\left( M\right)
\right) $ in the sequence $g_{\lambda _0^p}^{\perp }\subset g_{\lambda
_0^{p-1}}^{\perp }\subset \ldots \subset g_{\lambda _0^0}^{\perp }\subset g$
are subject to the  simple relation: 
$$\dim \left( N^{+}\left( so\left(
M\right) \right) \right) -\dim \left( N^{+}\left( so\left( M-4\right)
\right) \right) =2\left( \dim d_{so\left( M-4\right) }^v+1\right). 
$$ 
Taking
this into account we see that the chains (\ref{fullchain}) are full in the
sense that for $p=p^{\max }=\left[ M/4\right] +\left[ \left( M+1\right)
/4\right] $ their carrier spaces contain all the generators of $N^{+}\left(
so\left( M\right) \right) $. When $M$ is even-even or odd the total number
of Jordanian twists in a maximal full chain ${\cal F}_{{\cal G}_{0\prec
p\max }}$ is equal to the rank of $so\left( M\right) $. Thus in the latter
case the carrier subalgebra is equal to $B^{+}\left( so\left( M\right)
\right) $.

It was demonstrated in \cite{KS} how to construct new Yangians
using the explicit form of the twisting element. These new
Yangians are defined by the corresponding rational solution of the matrix
quantum Yang-Baxter equation (YBE). In particular, for the orthogonal
classical Lie algebras $so\left( M\right) $ one needs the twisting element $%
{\cal F}$ in the defining (vector) representation $d^v$ and the
auxiliary operators: the flip $P:v\otimes w\rightarrow w\otimes v$ ($P\in
{\rm Mat}\left( M\right) \otimes {\rm Mat}\left( M\right) $) and the
operator $K$, which is obtained from $P$ by transposing its first tensor
factor. The following expression gives the corresponding {\em deformed
rational solution} of the YBE:%
$$
ud^v\left( {\cal F}_{21}{\cal F}^{-1}\right) +P-\frac u{u-1+M/2}d^v\left(
{\cal F}_{21}\right) Kd^v\left( {\cal F}^{-1}\right)
$$
Here $u$ is a spectral parameter. In \cite{L} such deformed solutions were
obtained in the explicit form for the canonical 
chains ${\cal F}={\cal F}_{{\cal B}_{0\prec p}}$ .

All the calculations can be reproduced for the twisting 
elements ${\cal F}={\cal F}_{{\cal G}_{0\prec p}}$ of the full chains. 
This will lead to a new set
of so called {\em deformed Yangians} \cite{KST}.

\bigskip
{\small This work was partially supported by the Russian Foundation for
Basic Research under the grant 00-01-00500 (VDL) and 98-01-00310 (PPK).}

\end{document}